\documentclass[12pt,a4paper]{article}
\usepackage[cp1251]{inputenc} 
\usepackage[russian]{babel}
\usepackage{amsfonts}

\newcommand{\B}{$\hfill\Box$}
\newcommand{\al}{\alpha}

\newcommand{\la}{\lambda}

\newcommand{\vv}{\varphi}
\newcommand{\iy}{\infty}

\begin{document}

\begin{center}
{\large\bf
Inverse problems of recovering first-order integro-differential operators
from spectra}\\[0.2cm]
{\bf N.P.\,Bondarenko, V.A.\,Yurko} \\[0.3cm]
\end{center}

\thispagestyle{empty}

\noindent {\bf Abstract.} Inverse spectral problems are studied for
first-order integro-differential operators on a finite interval.
These problems consist in recovering some components of the kernel 
from one or multiple spectra.
Uniqueness theorems are proved for this class of inverse problems.



\begin{center}
{\bf 1. Introduction}
\end{center}

This paper deals with inverse spectral problems for first-order integro-diffferential operators in the form
$$
\ell y:=iy'(x)+\int_0^x M(x,t)y(t)\,dt,
$$
where $M(x,t)$ is a continuous function, called the kernel. Inverse problems of spectral analysis consist
in recovering operators from their spectral characteristics. 
For ordinary differential operators, such problems have been studied fairly completely (see [1-6] and the references
therein).

For integro-differential and other classes of nonlocal operators, inverse
problems are more difficult for investigation. The main classical methods
(transformation operator method and the method of spectral mappings)
either are not applicable to them or require essential modifications.
Therefore the general inverse problem theory for such operators has not been created yet.
At the same time, nonlocal and, in particular, integro-differential operators
are being actively studied, because they have many applications in physics, mechanics, biology, etc. (see, e.g., [7]).

Some fragmentary results of inverse problem theory for the first-order integro-differential operators
were obtained in [8,9] and other papers. However, those works deal with the kernels $M(x, t)$ of a very special structure,
while in the present paper, we consider a rather general case. We also mention the studies
[10-16], concerning inverse problems for second-order and higher-order integro-differential operators.

In this paper, we study inverse problems, that consist in recovering the kernel $M(x, t)$ from one or multiple spectra,
and prove the corresponding uniqueness theorems. Let us proceed with formulations of the main results.

Consider the boundary value problem $Q(M)$ for the integro-differential equation
$$
\ell y:=iy'(x)+\int_0^x M(x,t)y(t)\,dt=\la y(x),\quad x\in[0,\pi],    \eqno(1)
$$
with the condition $y(\pi) = 0$.

Together with $\ell$ we will consider the operator $\tilde\ell$ of
the same form but with a different kernel $\tilde M(x,t).$ We agree that
everywhere below if a certain symbol $\al$ denotes an object related
to $\ell,$ then $\tilde\al$ will denote the analogous object related
to $\tilde\ell.$

Fix $1 \le p \le \iy,$ and define the set of indices $I_p := \overline{1, p}$ for $p < \iy$ and $I_p := \mathbb N$ for $p = \iy.$
Let continuous functions
$R_j(x,t),\,0 \le t \le x \le \pi,\, j \in I_p$ be given, such that
$$
\int_0^x R_j(\pi-t,x-t)\,dt\ne 0,\quad 0<x\le\pi.
$$
Let
$$
M(x,t)=M_0(x,t)+\sum_{j=1}^p R_j(x,t)P_j(x-t),\quad
\tilde M(x,t)=M_0(x,t)+\sum_{j=1}^p R_j(x,t)\tilde P_j(x-t),
$$
where $M_0(x, t), P_j(x), \tilde P_j(x)$ are continuous functions for $0 \le t \le x \le \pi$, and for $p=\iy$
the series converge uniformly with respect to $x$ and $t.$ Let $\{\nu_{nk}\}$
be the eigenvalues (counted with multiplicities) of the boundary value problem
$Q_k:=Q(M_k)$ for $k \in I_p,$ and let $\{\tilde\nu_{nk}\}$ be the eigenvalues of the boundary value
problem $\tilde Q_k:=Q(\tilde M_k)$, where
$$
M_k(x,t)=M_0(x,t)+\sum_{j=1}^k R_j(x,t)P_j(x-t),\quad
\tilde M_k(x,t)=M_0(x,t)+\sum_{j=1}^k R_j(x,t)\tilde P_j(x-t).
$$

Our main result is the following uniqueness theorem.

\smallskip

{\bf Theorem 1. }{\it If $\nu_{nk}=\tilde\nu_{nk}$ for all $n$ and
$k \in I_p,$ then $M(x,t)\equiv\tilde M(x,t).$}

\smallskip

In particular, for $p = 1$ Theorem~1 takes the following form.

\smallskip

{\bf Theorem 2. }{\it Let the functions $M(x,t)$ and $\tilde M(x,t)$
have the form
$$
M(x,t)=M_0(x,t)+R(x,t)P(x-t),
\quad \tilde M(x,t)=M_0(x,t)+R(x,t)\tilde P(x-t),                   \eqno(2)
$$
where $M_0(x,t), R(x,t), P(x), \tilde P(x)$ are continuous functions and
$$
B(x):=\int_0^x R(\pi-t,x-t)\,dt \ne 0,\quad 0<x\le\pi,              \eqno(3)
$$
Let $\{ \nu_n \}$ and $\{ \tilde \nu_n \}$ be the eigenvalues (counted with multiplicities)
of the problems $Q(M)$ and $Q(\tilde M),$ respectively.
If $\nu_n=\tilde\nu_n$ for all $n,$
then $P(x)\equiv \tilde P(x)$ for $x\in[0,\pi].$}

\smallskip

Note that Theorem~1 can be easily derived from Theorem~2 by induction.
Therefore we will focus on the proof of Theorem~2.

The paper is organized as follows. In Section~2, we obtain the representation for the solution of Eq.~(1) 
in terms of the transformation operator. Note that transformation operators play an important role in spectral analysis
of both differential and integro-differential operators (see [1]). Section~3 is devoted to the proof of the uniqueness results.

\begin{center}
{\bf 2. Transformation operator}
\end{center}

Let $e(x,\la)$ be the solution of equation~(1), satisfying the initial condition $e(0,\la)=1.$

\smallskip
{\bf Lemma 1. }{\it The solution $e(x, \la)$ can be represented in the form
$$
e(x,\la)=\exp(-i\la x)+\int_0^x G(x,t)\exp(-i\la t)\,dt,            \eqno(4)
$$
where the function $G(x,t)$ is continuous for $0 \le t \le x,$ and}
$$
G(x,0)=0,\quad G(x,x)=i\int_0^x M(t,t)\,dt.
$$

{\it Proof. } One can easily show that the function $e(x,\la)$
satisfies the integral equation
$$
e(x,\la)=\exp(-i\la x)+i\int_0^x \exp(-i\la(x-t))\,dt
\int_0^t M(t,\tau)e(\tau,\la)\,d\tau.                              \eqno(5)
$$
Let us solve equation~(5) by the method of successive approximations:
$$
e(x,\la)=\sum_{n=0}^{\iy} e_n(x,\la),                            \eqno(6)
$$
where
$$
e_0(x,\la)=\exp(-i\la x),                                          \eqno(7)
$$
$$
e_{n+1}(x,\la)=i\int_0^x \exp(-i\la(x-t))\,dt
\int_0^t M(t,\tau)e_n(\tau,\la)\,d\tau,\quad k\ge 0.              \eqno(8)
$$
We will show by induction, that
$$
e_n(x,\la)=\int_0^x G_n(x,t)\exp(-i\la t)\,dt,\quad n\ge 1,       \eqno(9)
$$
where $G_n(x,t)$ are continuous functions, and $G_n(x,0)=0.$

Indeed, for $n=1$ we have
$$
e_1(x,\la)=i\int_0^x\,dt \int_0^t M(t,\tau)\exp(-i\la(x-t+\tau))\,d\tau
$$
$$
=i\int_0^x\,dt \int_{x-t}^x M(t,s+t-x)\exp(-i\la s)\,ds.
$$
Interchanging the order of integration, we obtain
$$
e_1(x,\la)=i\int_0^x \exp(-i\la s)\,ds \int_{x-s}^x M(t,s+t-x)\,dt
=\int_0^x G_1(x,s)\exp(-i\la s)\,ds,
$$
and consequently (9) is proved for $n=1$, where
$$
G_1(x,t)=i\int_{x-t}^x M(s,t+s-x)\,ds.                           \eqno(10)
$$
Suppose now that (9) is valid for a certain $n\ge 1.$
Then, substituting (9) into (8), we get
$$
e_{n+1}(x,\la)=i\int_0^x dt
\int_0^t M(t,\tau)\,d\tau\int_0^\tau G_n(\tau,\xi)\exp(-i\la\xi)\,d\xi
$$
$$
=i\int_0^x \exp(-i\la(x-t))\,dt \int_0^t M(t,\tau)\,d\tau
\int_{x-t}^{x-t+\tau} G_n(\tau,s+t-x)\exp(-i\la s)\,ds.
$$
Interchanging the order of integration, we obtain
$$
e_{n+1}(x,\la)=\int_0^x G_{n+1}(x,t)\exp(-i\la t)\,dt,
$$
where
$$
G_{n+1}(x,t)=i\int_{x-t}^x\,ds
\int_{t+s-x}^s M(s,\tau)G_n(\tau,t+s-x)\,d\tau.                 \eqno(11)
$$
The relation (10) and (11) imply that $G_n(x,0)=0$ for $n\ge 1.$
Substituting (9) into (6) and using (7), we arrive at (4), where
$$
G(x,t)=\sum_{n=1}^{\iy} G_n(x,t).                               \eqno(12)
$$
In view of (10) and (11), the series (12) converges
absolutely and uniformly for $0\le t\le x\le\pi,$ and the function
$G(x,t)$ is continuous. Moreover, $G(x,0)=0$ and
$$
G(x,x)=i\int_0^x M(t,t)\,dt.
$$
Lemma 1 is proved.\B

\begin{center}
{\bf 3. Uniqueness}
\end{center}

In this section, we prove the uniqueness Theorem~2. Theorem~1 follows from Theorem~2 by induction. 

The eigenvalues $\{\nu_n\}$ of the boundary value
problem $Q(M)$ coincide with the zeros of the entire function
$\Delta(\la):=e(\pi,\la).$ Using (4) and Hadamard's factorization
theorem, we obtain that the specification of the zeros $\{\nu_n\}$
uniquely determines the function $\Delta(\la).$
Under the assumptions of the theorem this yields
$$
e(\pi,\la)\equiv\tilde e(\pi,\la).
$$

Let the function $\psi(x,\la)$ be the solution of the equation
$$
\ell^{*}\psi:=-i\psi'(x,\la)+\int_x^{\pi} M(t,x)\psi(t,\la)\,dt=\la \psi(x,\la)  \eqno(13)
$$
under the condition $\psi(\pi,\la)=1.$ We multiply the relation
$\tilde\ell\tilde e(x,\la)=\la\tilde e(x,\la)$
by $\psi(x,\la)$, then subtract the relation (13), multiplied by
$\tilde e(x,\la)$, and integrate with respect to $x:$
$$
\int_0^\pi \psi(x,\la)\,dx \int_0^x (M(x,t)-\tilde M(x,t))
\tilde e(t,\la)\,dt=i(\tilde e(\pi,\la)-\psi(0,\la)).                \eqno(14)
$$
In particular, we get $\vv(\pi,\la)\equiv \psi(0,\la).$
We note that the function $w(x,\la):=\psi(\pi-x,\la)$ satisfies
the relations
$$
iw'(x,\la)+\int_0^x M(\pi-t,\pi-x)w(t,\la)\,dt=\la w(x,\la),\;w(0,\la)=1.
$$
Using (14), we derive
$$
\int_0^\pi \psi(x,\la)\,dx
\int_0^x (M(x,t)-\tilde M(x,t))\tilde e(t,\la)\,dt\equiv 0.       \eqno(15)
$$
The relations (2) and (15) imply that
$$
\int_0^\pi \psi(x,\la)\,dx
\int_0^x R(x,t)(P(x-t)-\tilde P(x-t))\tilde e(t,\la)\,dt\equiv 0. \eqno(16)
$$
The left-hand side of~(16) can be transformed as follows:
$$
\int_0^\pi \psi(x,\la)\,dx
\int_0^x R(x,t)(P(x-t)-\tilde P(x-t))\tilde e(t,\la)\,dt
$$
$$
=\int_0^\pi \psi(x,\la)\,dx
\int_0^x R(x,x-t)(P(t)-\tilde P(t))\tilde e(x-t,\la)\,dt
$$
$$
=\int_0^\pi (P(t)-\tilde P(t))\,dt \int_t^\pi
R(x,x-t)\psi(x,\la)\tilde e(x-t,\la)\,dx
$$
$$
=\int_0^\pi (P(\pi-x)-\tilde P(\pi-x))\,dx \int_0^x
R(\pi-t,x-t)w(t,\la)\tilde e(x-t,\la)\,dt.
$$
Consequently, the relation (16) implies
$$
\int_0^\pi (P(\pi-x)-\tilde P(\pi-x))z(x,\la)\,dx\equiv 0,      \eqno(17)
$$
where
$$
z(x,\la):=\int_0^x R(\pi-t,x-t)w(t,\la)\tilde e(x-t,\la)\,dt.  \eqno(18)
$$
In view of of Lemma~1, we have
$$
w(x,\la)=\exp(-i\la x)+\int_0^x K_1(x,t)\exp(-i\la t)\,dt,
$$
$$
\tilde e(x,\la)=\exp(-i\la x)+\int_0^x K_2(x,t)\exp(-i\la t)\,dt,
$$
where $K_j(x,t)$ are continuous functions. Substituting the latter relations into (18), we get
$$
z(x,\la)=B(x)\exp(-i\la x)
+\int_0^x R(\pi-t,x-t)\,dt\int_0^{t} K_1(t,\tau)\exp(-i\la(x-t+\tau))\,d\tau
$$
$$
+\int_0^x R(\pi-t,x-t)\,dt\int_0^{x-t} K_2(x-t,\xi)\exp(-i\la(t+\xi))\,d\xi
$$
$$
+\int_0^x R(\pi-t,x-t)\,dt\int_0^{t} K_1(t,\tau)\exp(-i\la\tau)\,d\tau
\int_0^{x-t} K_2(x-t,\xi)\exp(-i\la\xi)\,d\xi,
$$
and consequently,
$$
z(x,\la)=B(x)\exp(-i\la x)+\int_0^x K(x,t)\exp(-i\la t)\,dt,
$$
where $B(x)$ is defined in (3) and $K(x,t)$ is a continuous function.
Substituting the above expression into (17), we obtain
$$
\int_0^\pi (P(\pi-x)-\tilde P(\pi-x))
\Big(B(x)\exp(-i\la x)+\int_0^x K(x,t)\exp(-i\la t)\,dt\Big)\,dx\equiv 0.
$$
We rewrite the latter relation in the form
$$
\int_0^\pi \exp(-i\la x)\Big(B(x)(P(\pi-x)-\tilde P(\pi-x))+
\int_x^\pi K(t,x)(P(\pi-t)-\tilde P(\pi-t)) \Big)\,dx\equiv 0.
$$
Consequently, we get
$$
B(x)(P(\pi-x)-\tilde P(\pi-x))+
\int_x^\pi K(t,x)(P(\pi-t)-\tilde P(\pi-t))\equiv 0.
$$
Hence, in view of (3), we have $P(x)\equiv\tilde P(x).$
Theorem 2 is proved.

\bigskip
{\bf Acknowledgment.} This work was supported by Grant 17-11-01193
of the Russian Science Foundation.

\begin{center}
{\bf REFERENCES}
\end{center}
\begin{enumerate}
\item[{[1]}] Marchenko V.A., Sturm-Liouville operators and their applications.
     Naukova Dumka, Kiev, 1977;  English  transl., Birkh\"auser, 1986.
\item[{[2]}] Levitan B.M., Inverse Sturm-Liouville problems. Nauka,
     Moscow, 1984; English transl., VNU Sci.Press, Utrecht, 1987.
\item[{[3]}] Freiling G. and Yurko V.A., Inverse Sturm-Liouville Problems
     and their Applications. NOVA Science Publishers, New York, 2001.
\item[{[4]}] Beals R., Deift P. and Tomei C.,  Direct and Inverse Scattering
     on the Line, Math. Surveys and Monographs, v.28. Amer. Math. Soc.
     Providence: RI, 1988.
\item[{[5]}] Yurko V.A. Method of Spectral Mappings in the Inverse Problem
     Theory. Inverse and Ill-posed Problems Series. VSP, Utrecht, 2002.
\item[{[6]}] Yurko V.A. Inverse Spectral Problems for Differential
     Operators and their Applications. Gordon and Breach, Amsterdam, 2000.
\item[{[7]}] Lakshmikantham V. and Rama Mohana Rao M. Theory of
     integro-differential equations. Stability and Control: Theory and
     Applications, vol.1, Gordon and Breach, Singapure, 1995.
\item[{[8]}] Yurko V.A. Inverse spectral problems for first order integro-differential operators,
             Boundary Value Problems (2017) 2017:98.
\item[{[9]}] Buterin S.A. On an inverse spectral problem for first-order integro-differential
operators with discontinuities, Applied Mathematical Letters 78 (2018), 65--71.
\item[{[10]}] Yurko V.A. An inverse problem for integro-differential
     operators. Matem. zametki, 50, no.5 (1991), 134-146 (Russian);
     English transl. in Math. Notes, 50, no.5-6 (1991), 1188-1197.
\item[{[11]}] Yurko V. Inverse problems for arbitrary order integral and integro-differential operators, Results Math. (2018), 73:72.
\item[{[12]}] Kuryshova Yu. An inverse spectral problem for differential
     operators with integral delay. Tamkang J. Math. 42, no.3 (2011), 295-303.
\item[{[13]}] Buterin S.A. On the reconstruction of a convolution perturbation
     of the Sturm-Liouville operator from the spectrum, Diff. Uravn. 46 (2010),
     146--149 (Russian); English transl. in Diff. Eqns. 46 (2010), 150--154.
\item[{[14]}] Buterin S.A. and Vasiliev S.V. On uniqueness of recovering the convolution integro-differential operator from the spectrum
of its non-smooth one-dimensional perturbation, Boundary Value Problems (2018), 2018:15. 
\item[{[15]}] Zolotarev V.A. Inverse spectral problem for the operators with non-local potential, Mathematische Nachrichten (2018), 1--21,
DOI: https://doi.org/10.1002/mana.201700029.
\item[{[16]}] Wang Y., Wei G. The uniqueness for Sturm–Liouville problems with aftereffect, Acta Math. Sci. 32A:6 (2012), 1171--1178.

\end{enumerate}

\begin{tabular}{ll}
Name:             &   Bondarenko, Natalia \\
Place of work 1:  &   Department of Applied Mathematics and Physics, \\ 
{}                &   Samara National Research University \\
{}                &   Moskovskoye Shosse 34, Samara 443086, Russia \\
Place of work 2:  &   Department of Mechanics and Mathematics, Saratov State University\\
{}                &   Astrakhanskaya 83, Saratov 410012, Russia \\
E-mail:           &   bondarenkonp@info.sgu.ru
\end{tabular}

\bigskip

\begin{tabular}{ll}
Name:             &   Yurko, Vjacheslav \\
Place of work:    &   Department of Mechanics and Mathematics, Saratov State University\\
{}                &   Astrakhanskaya 83, Saratov 410012, Russia \\
E-mail:           &   yurkova@info.sgu.ru
\end{tabular}

\end{document}